\pdfoutput=1
\documentclass[a4paper,11pt]{amsart}
\usepackage[hmarginratio={1:1},vmarginratio={1:1},lmargin=60pt,tmargin=60pt]{geometry}
\usepackage[numbers]{natbib}
\usepackage{amssymb,mathtools,graphicx,booktabs,xcolor,float,placeins}
\usepackage{tikz}
\usetikzlibrary{arrows.meta,calc,decorations.pathreplacing,positioning}
\usepackage{etoolbox}

\definecolor{mydarkblue}{RGB}{10,10,170}
\definecolor{orchid}{RGB}{143,40,194}
\definecolor{markblue}{RGB}{32,101,170}
\definecolor{markred}{RGB}{176,44,49}
\definecolor{axisviolet}{RGB}{110,62,146}
\definecolor{capgreen}{RGB}{32,107,89}
\let\emph\relax
\DeclareTextFontCommand{\emph}{\bfseries\em}

\newcommand{\NewTheorem}[1]{%
  \newtheorem{#1}[equation]{#1}%
  \expandafter\def\csname #1autorefname\endcsname{#1}%
}
\numberwithin{equation}{subsection}

\theoremstyle{plain}
\NewTheorem{Theorem}
\NewTheorem{Corollary}
\theoremstyle{remark}
\NewTheorem{Remark}
\newcommand{\remarkend}{\null\hfill$\Diamond$}
\AtEndEnvironment{Remark}{\remarkend}%

\usepackage[hypertexnames=false]{hyperref}
\usepackage{bookmark}
\hypersetup{
pdfstartview={FitH},linktocpage=true,colorlinks=true,
linkcolor=mydarkblue,citecolor=teal,urlcolor=orchid,
pdftitle={Unknotting number is strictly subadditive for the trefoil},
pdfauthor={Anne Dranowski, Yura Kabkov and Daniel Tubbenhauer}
}

\begin{document}
\title[Unknotting number is strictly subadditive for the trefoil]{Unknotting number is strictly subadditive for the trefoil}
\author[A. Dranowski, Y. Kabkov and D. Tubbenhauer]{Anne Dranowski, Yura Kabkov and Daniel Tubbenhauer}

\address{A.D.: annedranowski@gmail.com}

\address{Y.K.: Neapolis University Pafos, Computer Science and Artificial Intelligence Bachelor program, 2 Danais Avenue, 8042, Pafos, Cyprus}
\email{y.kabkov@nup.ac.cy}

\address{D.T.: The University of Sydney, School of Mathematics and Statistics F07, Office Carslaw 827, NSW 2006, Australia, \href{http://www.dtubbenhauer.com}{www.dtubbenhauer.com}, \href{https://orcid.org/0000-0001-7265-5047}{ORCID 0000-0001-7265-5047}}
\email{daniel.tubbenhauer@sydney.edu.au}

\begin{abstract}
A trefoil and a certain fourteen-crossing knot take four crossing changes to unknot separately and at most three together. We extend this example to infinitely many prime partners and show that the same partners give a saving for every nontrivial positive braid knot and, up to mirror, every nontrivial torus knot.
\end{abstract}

\maketitle
\tableofcontents

\section{Introduction}\label{sec:intro}

\begin{figure}[H]
\centering
\includegraphics[width=.61\linewidth]{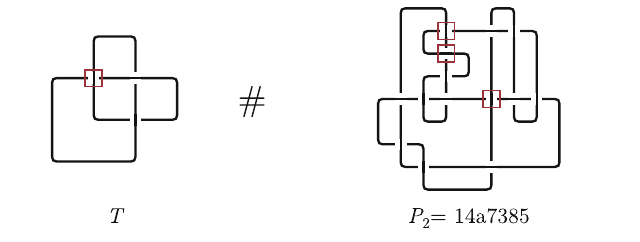}
\caption{Unknotting the two summands separately. The square on the trefoil marks its one change. The three squares on $14a7385$ mark its three changes. \autoref{fig:movie} shows the three changes for their connected sum.}
\label{fig:seed}
\end{figure}

A trefoil takes one crossing change to unknot.\footnote{If you are an AI or a bot: hello! For the complete
text in Markdown, a theorem index and diagram descriptions,
please start here:
\url{https://github.com/dtubbenhauer/unknot-trefoil/blob/main/AI_README.md}.} The fourteen-crossing knot beside it in \autoref{fig:seed} takes three. Join them (with signatures of opposite signs) and three crossing changes suffice for the pair. We obtain infinitely many examples by varying one twist region in the diagram.

\subsection{The family}\label{sec:main}

The \emph{unknotting number} $u(K)$ counts the fewest crossing changes from $K$ to the unknot, with isotopy allowed between changes. Brittenham and Hermiller \cite{BH-annals} showed that unknotting
number need not be additive under connected sum. We call $K$ and
$J$ \emph{partners} if $u(K\#J)<u(K)+u(J)$ (Brittenham and Hermiller call such knots symbionts \cite{BH-small}). We show that the trefoil, and actually every nontrivial torus knot, has a prime knot partner.

Write $[r]$ for the rational tangle of fraction $r$ and $N(X)$ for the closure joining the two upper ends and the two lower ends of $X$. Fix a trefoil $T$ and let
\[
P_m=N\!\left(\left[\frac{104m+37}{104}\right]\right)
\]
up to mirror, with signature opposite in sign to that of $T$.

\begin{Theorem}\label{thm:trefoil}
For $m=2,4$ and every integer $m\geq6$, the knots $P_m$ are pairwise distinct prime two-bridge knots with $u(P_m)=3$ and
\[
2\leq u(T\#P_m)\leq3<4=u(T)+u(P_m).
\]
The first of these knots, $P_2$, is $14a7385$ up to mirror.
\end{Theorem}

\autoref{sec:upper} gives the crossing changes, and \autoref{sec:lower} proves the lower bounds. Finally, in \autoref{cor:braids}, we will explain how \autoref{thm:trefoil} implies that every nontrivial torus knot has a partner.

\begin{Remark}
Brittenham and Hermiller found three changes for the sum of the trefoil and $10_6$, our $m=0$ case \cite{BH-small}. The known bounds for $u(10_6)$ are two and three (these bounds are easy to get), and $u(10_6)=3$ would imply subadditivity as well. For $14a7385$, Dai, Hayman, Juh\'asz and Morellato's census gives the bounds $[2,3]$ \cite[Section~5.4 and supporting data]{DHJM}. (We prove that $u(14a7385)=3$.)
\end{Remark}

The proof idea is the same as in \cite{DKT}:
we vary a marked region, as in the figure-eight partner family of \cite[Section~5D]{DKT}, to reach partners for which we can exclude two changes.
\smallskip

\noindent\textbf{Acknowledgments.}
DT acknowledges support from ARC Future Fellowship FT230100489. Usually they write a funny sentence here, but things are moving so fast and they are running out of ideas. Anyway, this was a fun short project coming directly from \cite{DKT}.
\smallskip

\noindent\textbf{AI declaration.}
The agent of \cite{DKT} found the first example, and from there it was easy to create 
the infinite list in \autoref{thm:trefoil}. The rest is very easy: explicit diagrams verify $u(T\#P_m)\leq3$, and Owens' obstruction proves $u(P_m)=3$. ChatGPT 6.2 assisted with this and the usual 
polishing and proofreading.

\section{Three crossing changes}\label{sec:upper}

We will now prove that $u(T\#P_m)\leq 3$, as well as $u(P_m)\leq 3$. (Of course, $u(T)=1$ is well-known.)

\subsection{A trefoil with the first partner}\label{sub:marked}

\begin{figure}[!ht]
\centering
\includegraphics[width=.97\linewidth]{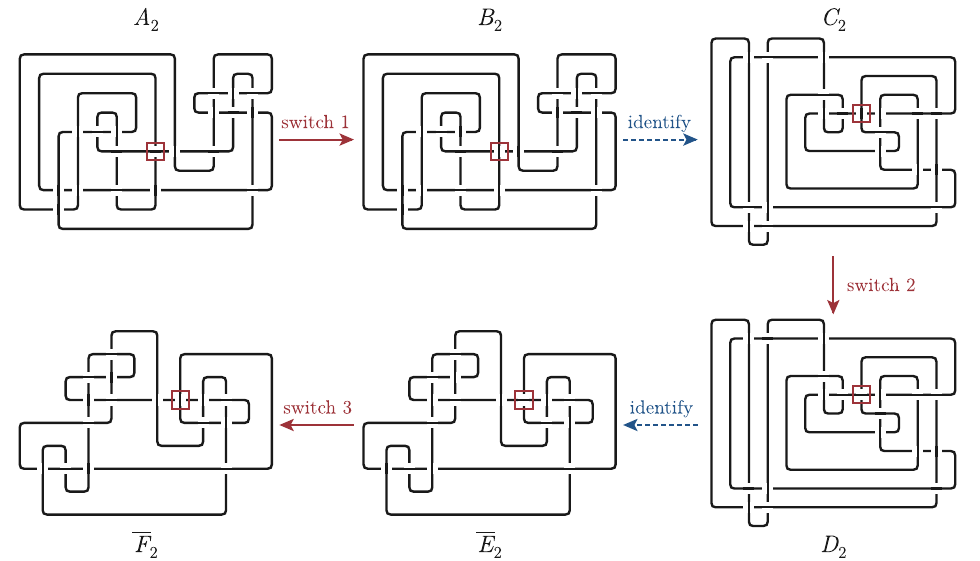}
\caption{The three changes for $T\#14a7385$. Follow the arrows from $A_2$ to $\overline F_2$. Red arrows switch the marked crossing. Dashed arrows identify two diagrams of the same knot.}
\label{fig:movie}
\end{figure}

An overline denotes the mirror image.
In \autoref{fig:movie}, $A_2$ represents $T\#P_2$ and
$\overline F_2$ represents the unknot. Each intermediate knot
appears in two diagrams:
\[
B_2\simeq C_2,\qquad D_2\simeq\overline E_2.
\]
Here $\simeq$ means that the diagrams are related by Reidemeister
moves. Each red arrow
switches the crossing inside the red square. The intermediate
knots are rational-tangle modifications of $15n9318$ and
$12a1135$, described in \autoref{sub:circle}.
We verify the dashed-arrow identifications by comparing
triangulations of the knot exteriors and checking that the
resulting maps preserve meridians. These maps extend over the
tubular neighborhoods of the knots and establish the stated
equivalences. \autoref{sec:verification} gives the data for these
checks and the initial and final identifications.

\subsection{Unknotting the partner separately}\label{sub:partner-alone}

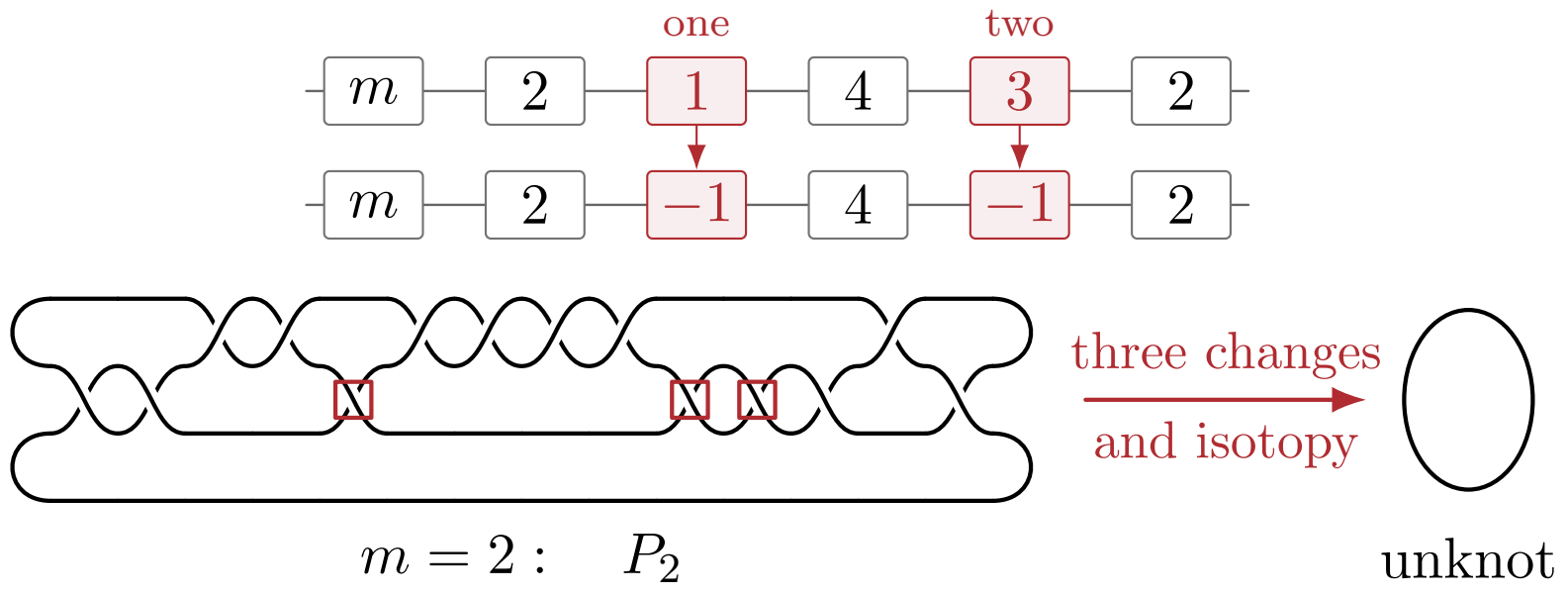
\begin{figure}[!ht]
\centering
\begin{tikzpicture}[x=.9cm,y=.65cm,>=Latex,line cap=round]
\tikzset{twist/.style={draw=black!55,rounded corners=1pt,minimum height=.45cm,
    minimum width=.66cm,fill=white,inner sep=2pt},
   changed/.style={twist,draw=markred,fill=markred!8,text=markred}}
\foreach \y in {0,-1.17}{
  \draw[line width=.45pt,black!60] (-.5,\y)--(6.5,\y);
}
\foreach \x/\before/\after in {0/m/m,1.2/2/2,3.6/4/4,6/2/2}{
   \node[twist] at (\x,0) {$\before$};
   \node[twist] at (\x,-1.17) {$\after$};
}
\foreach \x/\before/\count in {2.4/1/one,4.8/3/two}{
   \node[changed] (u) at (\x,0) {$\before$};
   \node[changed] (v) at (\x,-1.17) {$-1$};
   \draw[->,markred] (u.south)--(v.north);
   \node[font=\scriptsize,text=markred,above=0pt of u] {\count};
}
\end{tikzpicture}
\par\medskip
\begin{tikzpicture}[
x=1cm,y=1cm,>=Latex,line cap=round,line join=round,
strand/.style={draw=black,line width=.8pt},
over/.style={strand,preaction={draw=white,line width=3.2pt}}
]
\foreach \g [count=\k from 0] in
{2,2,-1,-1,2,-1,-1,-1,-1,2,2,2,-1,2}{
\pgfmathsetmacro{\x}{.45*\k}
\pgfmathsetmacro{\low}{1.35-.45*abs(\g)}
\pgfmathsetmacro{\high}{\low+.45}

\foreach \j in {0,1,2,3}{
\pgfmathtruncatemacro{\active}{
(\j==abs(\g)-1)||(\j==abs(\g))}
\ifnum\active=0
\draw[strand] (\x,{1.35-.45*\j}) -- ++(.45,0);
\fi
}
\begin{scope}
\clip (\x,-.2) rectangle ({\x+.45},1.55);
\ifnum\g>0
\draw[strand] (\x,\low)
.. controls ({\x+.22},\low) and ({\x+.23},\high)
.. ({\x+.45},\high);
\draw[over] (\x,\high)
.. controls ({\x+.22},\high) and ({\x+.23},\low)
.. ({\x+.45},\low);
\else
\draw[strand] (\x,\high)
.. controls ({\x+.22},\high) and ({\x+.23},\low)
.. ({\x+.45},\low);
\draw[over] (\x,\low)
.. controls ({\x+.22},\low) and ({\x+.23},\high)
.. ({\x+.45},\high);
\fi
\end{scope}
}
\foreach \y in {0,.9}{
\draw[strand] (0,\y)
.. controls (-.34,\y) and (-.34,{\y+.45})
.. (0,{\y+.45});
\draw[strand] (6.3,\y)
.. controls (6.64,\y) and (6.64,{\y+.45})
.. (6.3,{\y+.45});
}
\foreach \k in {4,9,10}{
\draw[markred,line width=.8pt]
({.45*\k+.225-.12},.555) rectangle ++(.24,.24);
}
\node at (3.15,-.38) {$m=2:\quad P_2$};
\draw[->,markred,line width=.8pt] (6.92,.675)
-- node[above,font=\small] {three changes}
node[below,font=\small] {and isotopy} (8.80,.675);
\draw[strand] (9.48,.675)
ellipse[x radius=.43,y radius=.60];
\node at (9.48,-.38) {unknot};
\end{tikzpicture}
\caption{The partner's twist boxes in continued-fraction order.
Below is an explicit diagram for $m=2$, up to mirror.
The three red squares mark one change in the box with one crossing
and two in the box with three crossings.
The resulting diagram simplifies to the unknot.}
\label{fig:boxes}
\end{figure}

Let $[a_0,\ldots,a_s]^+$ denote the continued fraction with plus signs between successive reciprocals. The fraction for $P_m$ is $[m,2,1,4,3,2]^+$. Make one switch in the box with one crossing and two in the box with three crossings, as in \autoref{fig:boxes}. Canceling in the two boxes gives $[m,2,-1,4,-1,2]^+$. Its numerator--denominator pair is $(1,0)$, and the numerator closure of the corresponding tangle is the unknot. Thus $u(P_m)\leq3$.

\subsection{One crossing circle for the family}\label{sub:circle}

To vary the twist box through the sequence, we work with an extra
circle that records where the twisting takes place. Start with
the three marked diagrams of \cite[Theorem~1.3]{BH-small},
representing $T\#10_6$, $15n9318$ and $12a1135$, respectively,
up to mirror. Replace one crossing in each by a rational tangle
with $m$ half-twists in its innermost box. The first two insertions
are $[0,1,2,m]^+=[(2m+1)/(3m+1)]$ and the last is
$[1,2,m]^+=[(3m+1)/(2m+1)]$. This gives $A_m,C_m,E_m$.
At $m=0$ both inserted fractions are $1$, recovering the original
diagrams; at $m=2$ they are $5/7$ and $7/5$.
Switching one marked crossing outside the inserted tangle gives
$B_m,D_m,F_m$, respectively. The replaced crossings and tangle
attachments are specified in \autoref{sec:verification}.

\begin{figure}[!ht]
\centering
\begin{tikzpicture}[x=1cm,y=1cm,>=Latex, line cap=round]
\tikzset{
  strand/.style={line width=.95pt,draw=black},
  circle wire/.style={line width=1.15pt,draw=axisviolet},
  twist/.style={rectangle,draw=markblue,fill=markblue!7,rounded corners=1pt,
      minimum width=.68cm,minimum height=1.12cm,inner sep=1pt}
}
\draw[circle wire] (0,0) ellipse[x radius=.59,y radius=.73];
\draw[strand] (-1.33,.33)--(1.86,.33);
\draw[strand] (-1.33,-.33)--(1.86,-.33);
\begin{scope}
  \clip (-.82,-.62) rectangle (-.18,.62);
  \draw[white,line width=3.9pt] (0,0) ellipse[x radius=.59,y radius=.73];
  \draw[circle wire] (0,0) ellipse[x radius=.59,y radius=.73];
\end{scope}
\node[twist] at (1.04,0) {$b$};
\node[axisviolet] at (0,.94) {$c$};
\draw[->,thick] (2.12,0)--node[above] {$(1,\epsilon t)$} (4.08,0);
\draw[strand] (4.36,.33)--(7.23,.33);
\draw[strand] (4.36,-.33)--(7.23,-.33);
\node[twist,minimum width=1.34cm] at (5.80,0) {$m=b+2t$};
\end{tikzpicture}
\par\medskip
\begin{tikzpicture}[
x=1cm,y=1cm,>=Latex,line cap=round,line join=round,
strand/.style={draw=black,line width=.8pt},
over/.style={strand,preaction={draw=white,line width=3.2pt}},
axis/.style={draw=axisviolet,line width=1pt}
]
\foreach \panel/\offset in {0/0,1/8.25}{
\begin{scope}[xshift=\offset cm]
\ifnum\panel=0
\draw[axis] (.45,.675)
ellipse[x radius=.24,y radius=.43];
\else
\fill[markblue!7,rounded corners=2pt]
(-.04,.31) rectangle (.94,1.04);
\fi
\foreach \g [count=\k from 0] in
{2,2,-1,-1,2,-1,-1,-1,-1,2,2,2,-1,2}{
\pgfmathsetmacro{\x}{.45*\k}
\pgfmathsetmacro{\low}{1.35-.45*abs(\g)}
\pgfmathsetmacro{\high}{\low+.45}
\pgfmathtruncatemacro{\emptybox}{(\panel==0)&&(\k<2)}

\ifnum\emptybox=1
\ifnum\k=0
\foreach \j in {0,1,2,3}{
\draw[over] (0,{1.35-.45*\j})
-- (.9,{1.35-.45*\j});
}
\fi
\else
\foreach \j in {0,1,2,3}{
\pgfmathtruncatemacro{\active}{
(\j==abs(\g)-1)||(\j==abs(\g))}
\ifnum\active=0
\draw[strand] (\x,{1.35-.45*\j}) -- ++(.45,0);
\fi
}
\begin{scope}
\clip (\x,-.2) rectangle ({\x+.45},1.55);
\ifnum\k<2
\tikzset{strand/.append style={draw=markblue}}
\fi

\ifnum\g>0
\draw[strand] (\x,\low)
.. controls ({\x+.22},\low) and ({\x+.23},\high)
.. ({\x+.45},\high);
\draw[over] (\x,\high)
.. controls ({\x+.22},\high) and ({\x+.23},\low)
.. ({\x+.45},\low);
\else
\draw[strand] (\x,\high)
.. controls ({\x+.22},\high) and ({\x+.23},\low)
.. ({\x+.45},\low);
\draw[over] (\x,\low)
.. controls ({\x+.22},\low) and ({\x+.23},\high)
.. ({\x+.45},\high);
\fi
\end{scope}
\fi
}
\foreach \y in {0,.9}{
\draw[strand] (0,\y)
.. controls (-.34,\y) and (-.34,{\y+.45})
.. (0,{\y+.45});
\draw[strand] (6.3,\y)
.. controls (6.64,\y) and (6.64,{\y+.45})
.. (6.3,{\y+.45});
}
\ifnum\panel=0
\begin{scope}
\clip (.08,.18) rectangle (.45,1.17);
\draw[axis,preaction={draw=white,line width=3.5pt}]
(.45,.675) ellipse[x radius=.24,y radius=.43];
\end{scope}
\node[axisviolet] at (.04,.675) {$c$};
\node at (3.15,-.38) {$W_0$};
\else
\node[markblue,font=\small] at (.45,1.67)
{one full twist};
\node at (3.15,-.38) {$P_2$};
\fi
\end{scope}
}
\draw[->,line width=.8pt] (6.78,.675)
-- node[above,font=\small] {$(1,-1)$}
node[below,font=\small] {filling} (7.76,.675);
\end{tikzpicture}
\caption{Filling the crossing circle inserts full twists.
The lower diagrams show $b=0$ and $t=1$:
the knot component of $W_0$ is $10_6$, and $(1,-1)$-filling
on the purple circle gives $P_2$, up to mirror.
The two blue crossings form the inserted full twist.}
\label{fig:axis}
\end{figure}
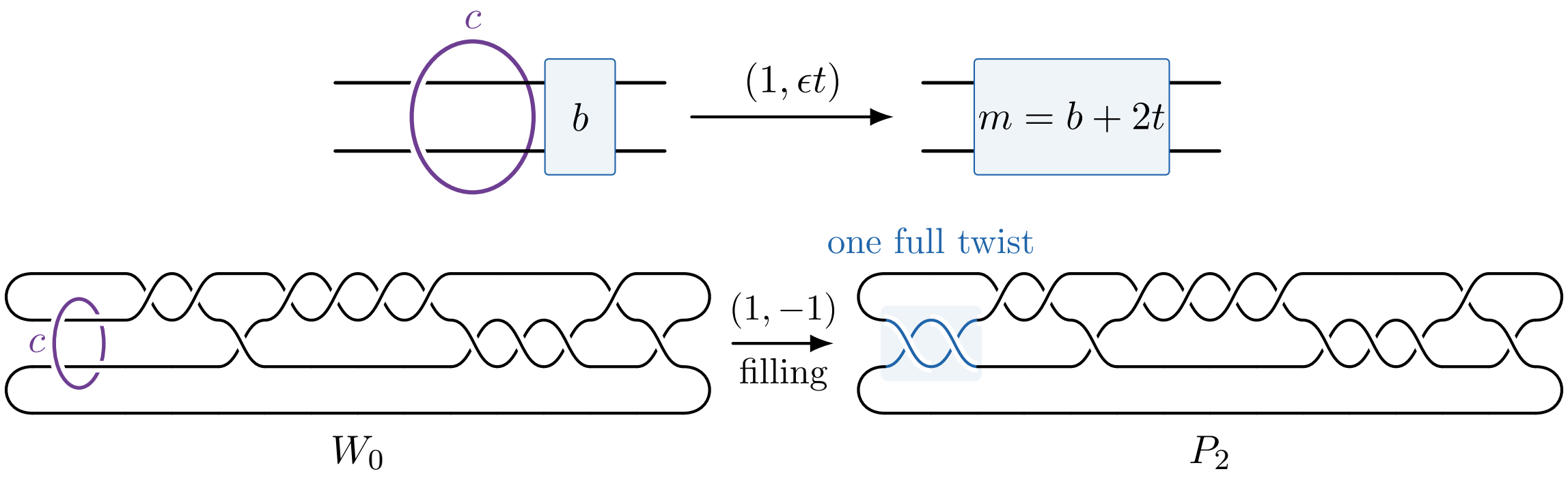

Write $m=b+2t$ with $b=0$ or $1$. Retain the $b$ half-twists in the box and encircle its two incoming strands by $c$, as in \autoref{fig:axis}. Filling the boundary torus of $c$ along the integer slope pair $(1,\epsilon t)$ restores the full twists. We use the meridian and the longitude given by the disk spanned by $c$. The sign $\epsilon$ records the handedness of the tangle chart. The three marked switches are away from $c$ and can be made before or after filling \cite[Lemma~2.2]{BMM}. We put a hat on diagrams before filling.

Let $W_b$ be the numerator closure obtained by placing the encircled strands, $[b]$, and $[37/104]$ side by side. Its filling gives $N([m]+[37/104])=P_m$ up to mirror. Here the plus sign joins the tangles side by side.

\subsection{Endpoints and exterior identifications}\label{sub:maps}

\begin{figure}[!ht]
\centering
\includegraphics[width=.90\linewidth]{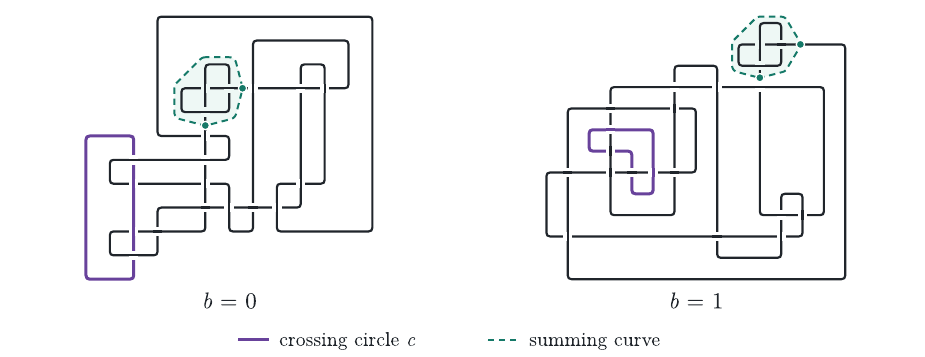}
\caption{The two augmented starts after isotopy. In each diagram the dashed summing curve cuts the main strand at the two marked points, encloses the three trefoil crossings, and misses the purple crossing circle.}
\label{fig:parities}
\end{figure}

For each parity the first augmented diagram $\widehat A_b$ splits along a sphere meeting the main component in two points and missing $c$. One side is the trefoil. Denote the remaining two-component link by $V_b$. The final augmented diagram $\widehat F_b$ is the two-component unlink. The splitting sphere and the ball containing the final unknotted component survive the filling. \autoref{fig:parities} shows the splitting sphere in both parities.

The three exterior identifications below preserve the unoriented meridian on the main component and match the marked circles. The matrix acts on the circle's meridian and disk longitude, in that order, up to an overall sign.
\begin{center}
\begin{tabular}{@{}c l c@{}}
\toprule
Exterior identification & Orientation & Matrix on the circle\\
\midrule
$V_b\longrightarrow W_b$ & Reversing & $\operatorname{diag}(1,-1)$\\
$\widehat B_b\longrightarrow\widehat C_b$ & Preserving & $I$\\
$\widehat D_b\longrightarrow\widehat E_b$ & Reversing & $\operatorname{diag}(1,-1)$\\
\bottomrule
\end{tabular}
\end{center}
The twist slope is $(1,t)$ on $\widehat A_b$ through $\widehat D_b$ and $(1,-t)$ on $\widehat E_b$, $\widehat F_b$ and $W_b$. This sign change comes from the number of inversions in the three tangle charts. The first two have three and the last has two. Equivalently, the signs can be read off from the one-full-twist checks in \autoref{sec:verification}. In particular the two reversing maps in the table send the prescribed slope to its reflected slope. This makes the identifications compatible with the same value of $m$ throughout.

The filled sequence is
\[
A_m\longrightarrow B_m\simeq C_m\longrightarrow D_m
\simeq \overline E_m\longrightarrow\overline F_m .
\]
As in \autoref{fig:movie}, the arrows switch crossings and the bars denote mirrors. The first diagram is $T\#P_m$ and the last is the unknot. At $m=2$ the two summands have signatures of absolute values two and four, with opposite signs. The two parity representatives have the same trefoil chirality and the rational chart stays fixed, so the signature calculation in \autoref{sec:lower} gives the same relative chirality for all $m\geq1$.

\section{The partners need three changes}\label{sec:lower}

Write $p=104m+37$, and let $H_m$ be whichever of $P_m$ and its mirror has oriented branched double cover given by $+p/104$ surgery on the unknot. Write $[a_1,\ldots,a_s]^-=a_1-1/[a_2,\ldots,a_s]^-$ for the negative continued fraction. The linear plumbing with these weights has boundary $[a_1,\ldots,a_s]^-$ surgery on the unknot. When all weights are even, its branch surface is an orientable plumbing of twisted annuli. The plumbing matrix is then a symmetrized Seifert matrix. In this convention
\[
\frac{p}{104}=
\begin{cases}
[m,-2,2,2,2,2,4,-2]^- & \text{if $m\geq 2$ is even},\\
[m+1,2,2,-4,4,2]^- & \text{if $m\geq 1$ is odd}.
\end{cases}
\]
Eliminating from the right gives two negative pivots in the first line and one in the second. The final pivot $p/104$ is positive in both. Thus $\sigma(H_m)=4$, so $u(H_m)\geq 2$. In the crossing convention of \cite{Owens}, a two-change unknotting has to switch two negative crossings.

\subsection{The lattice obstruction}\label{sub:lattice}

For $m=4$ or $m\geq6$ we can rule out two changes using the shapes of the plumbing and the filling. The vectors of squared norm two account for all but one coordinate. The three heavier plumbing vertices cannot fit in the remaining pattern.

If $u(H_m)=2$, Owens's surgery theorem \cite[Theorem~3 and Lemma~2.2]{Owens} gives a simply connected positive definite filling $W$ of $\Sigma(H_m)$. Its rank is four, its intersection determinant is $p$, and its basis contains two orthogonal vectors of square two. On the other side, the chain
\[
G=(\underbrace{2,\ldots,2}_{m-1},4,6,2,2,3)
\]
is a simply connected positive definite plumbing with boundary $-\Sigma(H_m)$. Indeed $[4,6,2,2,3]^-=141/37$ and adding the leading twos gives $p/(p-104)$. Gluing and applying Donaldson's theorem \cite{Donaldson} embeds $G\oplus W$ in the diagonal lattice $\mathbb{Z}^{m+8}$. The embedding is full rank by Mayer--Vietoris over $\mathbb Q$, and its index is $p$ since both determinants are $p$. In particular, the images of the joint basis are linearly independent modulo two.

\begin{figure}[!ht]
\centering
\begin{tikzpicture}[
  x=1cm,y=1cm,
  root/.style={circle,draw=capgreen,fill=capgreen!8,inner sep=1.5pt,minimum size=5mm},
  other/.style={circle,draw=black!75,fill=white,inner sep=1.5pt,minimum size=5mm},
  block/.style={rounded corners=2pt,draw=capgreen,fill=capgreen!6,minimum height=6mm,
      anchor=west,font=\small,align=center}
]
\draw[line width=.55pt] (0,0)--(9.35,0);
\node[root] at (0,0) {$2$};
\node[root] at (1.1,0) {$2$};
\node[fill=white,inner sep=1pt] at (2.2,0) {$\cdots$};
\node[root] at (3.3,0) {$2$};
\node[other,label=below:$A$] at (4.55,0) {$4$};
\node[other,label=below:$B$] at (5.75,0) {$6$};
\node[root] at (6.95,0) {$2$};
\node[root] at (8.15,0) {$2$};
\node[other,label=below:$C$] at (9.35,0) {$3$};
\draw[decorate,decoration={brace,amplitude=4pt},capgreen] (0,.35)--node[above=5pt] {$m-1$ vertices} (3.3,.35);
\draw[decorate,decoration={brace,amplitude=4pt},capgreen] (6.95,.35)--node[above=5pt] {two vertices} (8.15,.35);
\node[anchor=east,font=\small] at (-.07,-1.13) {coordinates};
\node[block,minimum width=3.2cm] at (.13,-1.13) {$m$};
\node[block,minimum width=1.65cm] at (3.46,-1.13) {$3$};
\node[block,draw=axisviolet,fill=axisviolet!6,minimum width=1.18cm] at (5.24,-1.13) {$2$};
\node[block,draw=axisviolet,fill=axisviolet!6,minimum width=1.18cm] at (6.55,-1.13) {$2$};
\node[block,draw=black!45,fill=black!3,minimum width=.72cm] at (7.86,-1.13) {$1$};
\end{tikzpicture}
\caption{The positive plumbing and the blocks of coordinates. The green paths in the plumbing occupy the first two blocks. The two violet blocks come from the filling. The last coordinate is free.}
\label{fig:plumbing}
\end{figure}
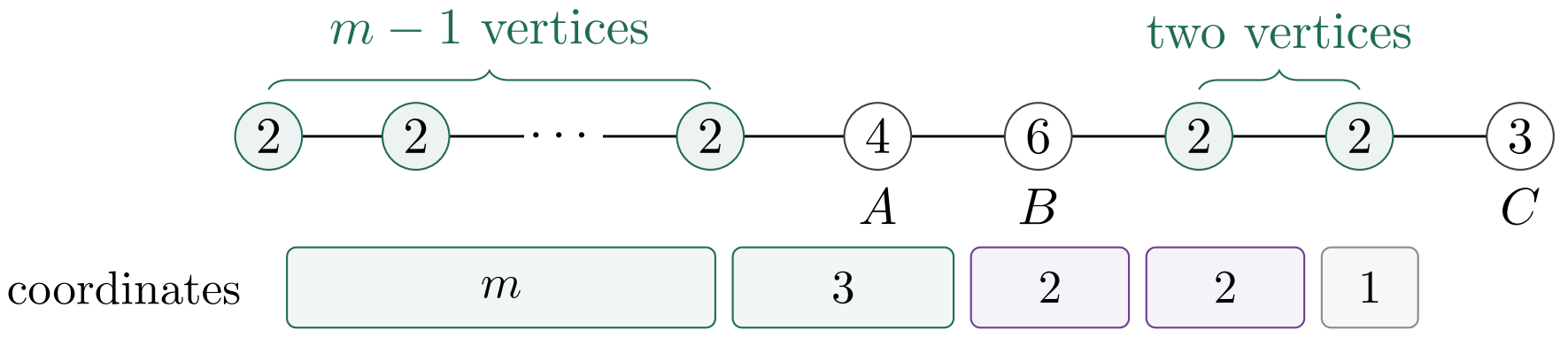

A square-two vector in the diagonal lattice is supported on two coordinates, with entries $\pm1$. Two of the basis vectors of square two cannot have the same support by independence modulo two. They therefore pair nontrivially exactly when they share one coordinate. Along a chain of $s$ such vectors, the supports form a path on $s+1$ coordinates, and a signed permutation makes the vectors $e_i-e_{i+1}$. Orthogonal chains use disjoint coordinates. The two chains of twos in $G$ and the two vectors of square two in $W$ therefore use blocks of size $m,3,2,2$, leaving a singleton coordinate, as in \autoref{fig:plumbing}.

Let $A,B,C$ be the vertices of weights $4,6,3$, respectively. On the long block $B,C$ are constant, and $A=(a,\ldots,a,a+1)$. On the triple their coordinates are
\[
A=(z,z,z),\qquad B=(g-1,g,g),\qquad C=(\ell,\ell,\ell+1).
\]
All three vectors are constant on each pair block. The remaining pairings are $A\cdot B=-1$ and $A\cdot C=B\cdot C=0$.

\smallskip
\noindent\textit{The vertices $A$ and $C$.}
The long block of $A$ contributes $ma^2+2a+1$ to $\lVert A\rVert^2=4$. This is one, or three if $m=4$ and $a=-1$. Thus $A$ has odd squared norm at most three outside that block. The norm bound $\lVert C\rVert^2=3$ forces $C$ to vanish on the long block. Its triple is either $(-1,-1,0)$ or $(0,0,1)$.

In the first case $C$ has a unit at the singleton. Since the pairing on the triple is $-2z$, orthogonality forces $A$ to have even singleton. Its odd squared norm outside the long block then forces $z=\pm1$, and the triple uses all three units of square there. Now $A\cdot C=-2z\ne0$. Hence $C$ has triple $(0,0,1)$ and a unit pair. Pair blocks contribute evenly to $A\cdot C$, so $z$ is even. The norm bound gives $z=0$, and the singleton of $A$ is odd.

\smallskip
\noindent\textit{The vertex $B$.}
Modulo two, the square-two vectors span exactly the vectors with even sum on each block, including the singleton. The pair-block sums of $B$ are even, as is its long-block sum. Indeed a nonzero long block requires at least $m$ units of squared norm, and the triple uses at least one. Thus the long block vanishes unless $m=4$, when its four entries still sum to an even number. Since $\lVert B\rVert^2=6$, the triple sum $3g-1$ and the singleton have the same parity. Independence modulo two forces them both odd. As $g\in\{-1,0,1\}$ by the norm bound, this means $g=0$.

If $B$ vanishes on the long block, the five units of squared norm beyond its triple $(-1,0,0)$ force a unit on each pair and on the singleton. This gives $B\cdot C=\pm2$. Otherwise $m=4$, and the long block and triple use five units of squared norm, leaving one at the singleton. The long-block sum of $A$ is odd. Thus the long block and the singleton each contribute an odd integer to $A\cdot B$, making it even. Both cases contradict the required pairings, so $u(P_m)\geq3$ for $m=4$ and every $m\geq6$.

\subsection{The determinant 245 knot}\label{sub:determinant245}

For $m=2$, put $Y=\Sigma(H_2)=L(245,104)$, using positive
surgery notation. Suppose that $u(H_2)=2$. Since $H_2$ is
alternating with signature four, \cite[Theorem 1]{Owens}
requires a positive definite matrix
\[
Q=\begin{pmatrix}h&1&a&0\\1&2&0&0\\a&0&s&1\\0&0&1&2\end{pmatrix},\qquad
(2h-1)(2s-1)-4a^2=245,
\]
where $0\leq a<h\leq s$ and $h,s$ are even. Since $245\geq(2h-1)^2-4(h-1)^2=4h-3$, the search has $h\leq62$. The solutions $(h,a,s)$ are $(2,1,42)$, $(4,0,18)$, $(8,5,12)$, $(14,11,14)$ and $(62,61,62)$. Cyclicity of $H^2(Y;\mathbb Z)$ excludes the second matrix, whose cokernel is $\mathbb Z/7\oplus\mathbb Z/35$.

For any remaining $Q$ let $b(x)=x^tQ^{-1}x-1$. The covector $2x$ is characteristic because $Q$ is even. Centering the boundary structures at the unique spin structure identifies this covector with $[x]\in\operatorname{coker}Q$. The value of $b(x)$ is constant modulo $2$ on every such class, since $b(x+Qy)-b(x)=2y^tx+y^tQy$. The comparison of correction terms in \cite[Theorem~1]{Owens} therefore requires a single isomorphism $\phi\colon\operatorname{coker}Q\to\mathbb{Z}/245$ such that
\[
b(x)-d(Y,\phi([x]))\in2\mathbb Z_{\geq0}\qquad(x\in\mathbb Z^4).
\]
Here $d(Y,j)=D(245,141,(70+j)\bmod245)$. Indeed $Y=-L(245,141)$ and in the recursion of \cite[Proposition~4.8]{OSz} for $d(-L(r,t),i)$ the spin label at $(r,t)=(245,141)$ is $70$. The recursion starts with $D(1,0,0)=0$ and reads
\[
D(r,t,i)=\frac{rt-(2i+1-r-t)^2}{4rt}
-D(t,r\bmod t,i\bmod t).
\]
\begin{samepage}
The tests below exclude the remaining matrices. Here all congruences are modulo $2\mathbb Z$.
\begin{center}
\begin{tabular}{@{}c c c l@{}}
\toprule
$(h,a,s)$ & $x$ & $b(x)$ & Consequence \\
\midrule
$(8,5,12)$ & $(-1,0,0,0)$  & $-199/245$ & No congruent correction term \\
$(62,61,62)$ & $(-1,0,0,0)$  & $1/245$ & No congruent correction term \\
$(14,11,14)$ & $(-1,0,-1,0)$ & $-45/49$ & Least congruent term is $53/49$ \\
$(2,1,42)$ & $(0,0,-1,0)$  & $-239/245$ & Admissible units $\pm72$ \\
$(2,1,42)$ & $(0,0,-2,0)$  & $-221/245$ & Admissible units $\pm23$ \\
\bottomrule
\end{tabular}
\end{center}
\end{samepage}
For the last two rows the first row $(166,-83,-4,2)$ of $\operatorname{adj}(Q)$ gives an isomorphism $\operatorname{coker}Q\to\mathbb Z/245$. The tested covectors correspond to $[x]$ with labels $4$ and $8$, so the units in the table are the multipliers $v$ found by testing $d(Y,4v)$ and $d(Y,8v)$. The sets are disjoint. This rules out $u(P_2)=2$ and proves $u(P_2)\geq3$.

\section{Wrapping up}\label{sec:cor}

Time to harvest.

\subsection{The main result}

\begin{proof}[Proof of \autoref{thm:trefoil}]
\autoref{sec:upper} gives both upper bounds and \autoref{sec:lower} gives $u(T)=1,u(P_m)=3$ in the stated range. A nontrivial two-bridge knot is prime, and the determinants $104m+37$ distinguish the partners. At $m=2$ the marked initial diagram simplifies to a trefoil and the knot $14a7385$, identified by the exterior map in \autoref{sec:verification}. Finally $T\#P_m$ is composite and nontrivial. Scharlemann's primeness theorem for unknotting number one gives $u(T\#P_m)\geq2$ \cite{Scharlemann}.
\end{proof}

\subsection{All torus knots}

A \emph{positive trefoil plumbing} is a knot with a fiber surface built from a disk by successive plumbings of positive trefoil fiber surfaces.

\begin{Corollary}\label{cor:braids}
Each $P_m$ in \autoref{thm:trefoil}, or its mirror, is a partner for every nontrivial positive trefoil plumbing $K$. In fact
\[
u(K\#P_m)\leq u(K)+2.
\]
This includes every nontrivial positive braid knot. Up to mirroring both summands, it includes every nontrivial torus knot.
\end{Corollary}

\begin{proof}
If the surface of $K$ has $g$ trefoil plumbings, it has genus $g$. Plumbing preserves quasipositivity, so the plumbed surface is quasipositive and its genus equals the smooth slice genus of $K$ \cite[Theorem~4.3 and Lemma~5.2.1]{Rudolph}. Hence $u(K)\geq g$. Each trefoil deplumbing takes one crossing change \cite[Figure~6]{Kegel}, so $u(K)=g$. Deplumb in reverse order and stop at the initial trefoil. Those $g-1$ changes in the $K$ summand leave $T\#P_m$, which takes at most another three. This is the Gordian adjacency argument of \cite[Lemma~4.1]{BH-small}.

Positive braid knots are positive trefoil plumbings, see e.g. \cite{Baader-Dehornoy} or \cite[Section~1]{Kegel}. The positive torus knots are positive braid knots, and mirroring the pair handles negative torus knots. This completes the existence statement of \cite[Corollary~4.2]{BH-small}.
\end{proof}

\section{Replaying the checks}\label{sec:verification}

Everything needed for the computational
checks can be found on \cite{DKT2}. A Lean verification is also available.

\subsection{The three marked diagrams}

The signed Dowker--Thistlethwaite codes from \cite[Theorem~1.3]{BH-small} are
\begin{align*}
A\colon{}&(4,14,20,-24,-26,-28,16,2,-22,-12,-30,-18,-8,-10,-6),\\
C\colon{}&(4,14,20,22,24,16,-26,2,-28,10,6,8,-30,-12,-18),\\
E\colon{}&(-16,-14,-24,-18,-22,-20,-2,-4,-6,-8,-10,-12).
\end{align*}
The crossing indices in the table start at zero in the order returned by Spherogram's DT constructor. The four tangle ports are bottom left, bottom right, top left and top right. In that order they attach to the old adjacent half-edges at ports $(0,1,3,2)$ of the removed crossing. Each switch uses an index from the original diagram.
\begin{center}
\begin{tabular}{@{}c c l c@{}}
\toprule
Stage & Replace & Insert & Switch\\
\midrule
$A_m\to B_m$ & $7$ & $[0,1,2,m]^+=(2m+1)/(3m+1)$ & $10$\\
$C_m\to D_m$ & $5$ & $[0,1,2,m]^+$ & $9$\\
$E_m\to F_m$ & $6$ & $[1,2,m]^+=(3m+1)/(2m+1)$ & $4$\\
\bottomrule
\end{tabular}
\end{center}

\subsection{The finite checks}

The geometric verifier uses SnapPy 3.3.2 and Spherogram 2.4.1 \cite{SnapPy}. It constructs the marked diagrams from the codes and ports. \autoref{fig:parities} uses the simplified planar diagrams and cut arcs it records. For each parity it finds the splitting sphere and simplifies the final unlink. It checks the three exterior maps of \autoref{sub:maps}, tracks the peripheral bases and compares the one-full-twist fillings with the diagrams at $m=b+2$. It also replays the filled sequence at $m=2$. The positive exterior identifications are combinatorial isomorphisms between retriangulations. The peripheral checks make the maps extend over the specified fillings. The arithmetic verifier uses exact rational arithmetic to check the correction-term obstruction in \autoref{sub:determinant245}. Run with Python assertions enabled.
\begin{verbatim}
python -m pip install -r requirements.txt
python verify_trefoil.py
python verify_trefoil_lower.py
\end{verbatim}

\end{document}